\documentclass{gtart}
\usepackage{amsmath,amssymb,amscd}

\input gtmonout
\volumenumber{2}
\volumeyear{1999}
\volumename{Proceedings of the Kirbyfest}
\pagenumbers{473}{488}
\papernumber{23}
\received{14 October 1999}\revised{18 November 1999}
\published{23 November 1999}

\theoremstyle{plain}
\newtheorem{theorem}{Theorem}[section]

\newtheorem{ther}{Theorem}

\theoremstyle{remark}
\newtheorem*{remark}{Remark}

\theoremstyle{definition}
\newtheorem{defn}[theorem]{Definition}

\def \R {\mathbf{R}}
\def \Z {\mathbf{Z}}

\def \D {\mathbf{D}}
\def\cp{\mathbf{CP}}
\def\cpbar{\overline{\mathbf{CP}}^{2}}
\def \HH{\mathcal{H}}
\def \WW{\mathcal{W}}

\def\LL{\mathcal{L}}
\def\BB{\mathcal{B}}
\def\DD{\mathcal{D}}

\def\MM{\mathcal{M}}
\def\VV{\mathcal{V}}

\def\diff{\mathit{Diff}}
\def\diffp{\mathit{Diff}^+}
\def\diffh{\mathit{Diff}_H}
\DeclareMathOperator{\SU}{SU}
\DeclareMathOperator{\SO}{SO}
\DeclareMathOperator{\met}{Met}

\def\SW{\ifmmode{\text{SW}}\else{$\text{SW}$}\fi}
\def\SWW{\ifmmode{\underline{\text{SW}}}\else{$\underline{\text{SW}}$}\fi}

\def\S{section }

\begin{document}

\title{A polynomial invariant of diffeomorphisms\\of 4--manifolds}
\asciititle{A polynomial invariant of diffeomorphisms of 4-manifolds}
\author{Daniel Ruberman}
\address{Department of Mathematics, Brandeis University\\
Waltham, MA 02254, USA}
\email{ruberman@binah.cc.brandeis.edu}
\begin{abstract}
We use a 1--parameter version of gauge theory to investigate the
topology of the diffeomorphism group of 4--manifolds.   A polynomial invariant,
analogous to the Donaldson polynomial, is defined, and is used to show that the
diffeomorphism group of certain simply-connected 4--manifolds has infinitely
generated
$\pi_0$.
\end{abstract}
\asciiabstract{We use a 1-parameter version of gauge theory to investigate the
topology of the diffeomorphism group of 4-manifolds.   A polynomial invariant,
analogous to the Donaldson polynomial, is defined, and is used to show that the
diffeomorphism group of certain simply-connected 4-manifolds has infinitely
generated \pi_0.}
\maketitle
\cl{\small\it Dedicated to Rob Kirby on the occasion of his $60^{th}$ birthday.}

\section{Introduction}
The issue of whether topological and smooth isotopy coincide for
diffeomorphisms of $4$--manifolds was recently resolved in the author's
paper~\cite{ruberman:isotopy}.  That work defined an invariant, roughly
analogous to the degree--$0$ part of the Donaldson invariant of a $4$--manifold,
which serves as an effective obstruction to smooth isotopy.  In the current
paper, we will extend the definition of the invariant to give a polynomial-type
invariant, which is analogous to the full Donaldson polynomial.  As an
application of the polynomial invariant, we will show that $\pi_0$ of the
diffeomorphism group of certain $4$--manifolds is infinitely generated.

It is worth stating this last result somewhat more precisely.   For any compact
$4$--manifold $X$, one can consider its (orientation-preserving) diffeomorphism
group $\diffp(X)$.  Taking the induced map on homology defines a homomorphism
from
$\diffp(X)$ to the automorphism group of the intersection form of $X$; in many
cases this map is a surjection.  Let us denote by $\diffh(X) \subset \diffp(X)$
the kernel of this map.
\begin{ther}\label{infgen} Let $Z_n$ $(n\ge 2)$ denote the connected sum
$$
\#_{k_n} \cp^2 \#_{l_n} \cpbar
$$
where $k_n = 2n$ and $l_n = 10n+1$.  Then there is a homomorphism
$$
\D\co \pi_0(\diffh(Z_n))  \to \R[[ H_2(Z_n)^* ]]
$$
with infinitely generated image.
\end{ther}
The numbers appearing in the definition of $Z_n$ are less obscure than might
appear at first glance; the manifold $Z_n$ is diffeomorphic to the elliptic
surface
$E(n)$, connected sum with $\cp^2$ and two copies of $\cpbar$.  As will become
evident in the proof, the conclusion that the image of $\D$ is infinitely
generated derives from the fact that
$E(n)$ supports infinitely many smooth structures which become diffeomorphic
upon connected-sum with $\cp^2$.\\
{\bf Acknowledgement}\qua  The author was partially supported by
NSF Grant 4-50645.

\section{Invariants of diffeomorphisms}
Let us start with a brief review of the definition of the $0$--degree invariant
discussed in~\cite{ruberman:isotopy}.   The conditions discussed below having
to do with orientations are used in defining the invariant as an element of
$\Z$, rather than merely modulo $2$. The data necessary for the definition are:
\begin{enumerate}
\item A smooth, simply-connected, oriented, homology-oriented $4$--manifold $Y$
with $b_+^2 > 2$.
\item An $\SO(3)$ bundle $P \to Y$ such that $w_2(P)\neq 0$, and with
$\dim(\MM(P)) = -1$.  (Here $\MM(P)$ is the moduli space of anti-self-dual
connections on $P$.)
\item An integral lift $c\in H^2(Y;\Z)$ of $w_2(P)$.
\item An orientation-preserving diffeomorphism $f$ of $Y$ such that $f^*(P)
\cong P$, and such that the quantity $\alpha(f)\beta(f) = 1$.
\end{enumerate}
The product $\alpha\beta \in \{\pm 1\}$  in the last item indicates, roughly,
whether $f$ preserves or reverses the orientation of the moduli space.   The
numbers $\alpha$ and $\beta$ are themselves defined as follows:
\begin{itemize}
\item
Composing a projection of $H^2(Y;\R)$ onto $H^2_+(Y)$ with $f^*$ defines an
isomorphism of $H^2_+(Y)$ with itself; the sign of the determinant (which is
independent of all choices) determines the spinor norm, $\alpha(f) \in \{\pm
1\}$.
\item
The condition that $f^*P \cong P$ implies that $f^*w_2 = w_2$, or in other
words that $f^*c -c$ is divisible by $2$ in
$H^2(Y)$.  One thereby can define $\beta(f) =
(-1)^{(\frac{f^*c-c}{2})^2}$.
\end{itemize}
Under these conditions, for a generic metric $g\in \met(Y)$, the moduli space
$\MM(P;g_0)$ (connections which are $g_0$--anti-self-dual) is empty.  If one
considers instead a generic path $g_t \in \met(Y)$ of metrics from $g_0$ to
$g_1=f^*g_0$, then one can construct the $1$--parameter
moduli space
$$\tilde\MM(P;\{g_t\}) = \bigcup_{t\in [0,1]} \MM(P;g_t).
$$
The count of points, with signs, in this $0$--dimensional moduli space defines
an invariant $D(f)$ (or $D_Y(f;P)$ if one needs to keep track of the
manifold and/or the bundle).

The independence of $D(f)$ from the choice of initial metric $g_0 \in \met(Y)$
and of the choice of generic path are proved using a
$2$--parameter moduli space
$$
\Tilde{\Tilde{\mathcal{M}}}  = \bigcup_{(s,t)\in I \times I} \MM(P,K_{s,t}).
$$
Here $K_{s,t}$ is a $2$--parameter family of metrics giving a homotopy from one
path of metrics $g_s=K_{s,0}$ to $k_s=K_{s,1}$.   The proof in each case uses
a choice of `boundary conditions' for the endpoints of the homotopy.  A
fundamental point is that the parameter space $\met(Y)$ is simply connected, so
that an arbitrary assignment of metrics on the boundary of the
$(s,t)$ square $I
\times I$ can be filled in smoothly.   So for instance, to verify independence
from the choice of path, use a $2$--parameter family in which the endpoints are
fixed:
$K_{0,t}=g_0$ and $K_{1,t}=g_1$.  To verify that the initial metrics  $g_0$ and
$k_0$ give the same value for $D(f)$, use an arbitrary path from $g_0$ to $k_0$
for $K_{0,t}$ with the proviso that the right endpoints $K_{1,t}$ are equal to
$f^*K_{0,t}$.

In both arguments, the principle used is that on the one hand, the
boundary of the
$2$--parameter moduli space $\Tilde{\Tilde{\mathcal{M}}}$, which is a compact
$1$--manifold, consists of algebraically $0$ points.   On the other hand, the
boundary is also the union of the $1$--parameter moduli spaces associated to the
four sides of the $(s,t)$ square.  In the first case, the right and left sides
of the square are fixed at generic metrics defining empty moduli spaces, so the
boundary is the difference between the invariant computed with the two
different paths on the top and bottom.  In the second case, one must account
for the additional part of the boundary, given by the difference between the
(algebraic) count of points on the left and right sides.  However, the choice
of boundary conditions takes care of this, because there is an isomorphism
between $\tilde\MM(P;K_{0,t})$ and
$\tilde\MM(P;f^*K_{0,t})$ and so the contributions from the two sides cancel.

\subsection{A polynomial invariant of diffeomorphisms}\label{polysec}
It is natural to try to extend $D(f)$ to a polynomial in $H_0(Y) \oplus
H_2(Y)$, by considering an $\SU(2)$ or $\SO(3)$ bundle $P$ for which the
ASD moduli space
$\MM(P)$ has positive odd dimension, and cutting down by divisors.  Recall that
in the construction of the usual Donaldson polynomial, the divisor
associated to
a surface $\Sigma$ in $Y$ is defined in several steps.  One first considers the
space of irreducible connections on $\Sigma$, together with the restriction map
$r_\Sigma\co  \BB^{*,\Sigma}(Y) \to \BB^*(\Sigma)$. Here $\BB^{*,\Sigma}(Y) $
consists of connections whose restriction to $\Sigma$ is irreducible.  An
important remark (cf~\cite[\S 9.2.3]{donaldson-kronheimer}) is that for a
generic surface
$\Sigma \subset Y$, the moduli space $\MM(P)$ is contained in
$\BB^{*,\Sigma}(Y)$.  There is a natural line bundle
$\LL \to
\BB^*(\Sigma)$, and one chooses a section, which is then pulled back to
$\BB^{*,\Sigma}(Y)$.  If these constructions are done with some care, then the
zero-set of the pulled-back section defines a divisor $V_\Sigma$.   Since we
are only concerned with intersections of $V_\Sigma$ with $\MM(P)$, we will
follow the standard notational abuse and drop the superscript $\Sigma$. A
similar construction gives a codimension--$4$ submanifold
$V_x$ of
$\BB^*(Y)$ which represents the dual of $\mu$ of a point $x\in Y$.

Now these constructions depend on a number of choices, e.g. the specific
representative of the homology class $[\Sigma]$, and the choice of section of
$\LL_\Sigma$.   If the `space' of possible choices were simply-connected, then
one could incorporate them into the parameter space
$\Pi$, and proceed precisely as in the definitions in Section 2
of~\cite{ruberman:isotopy}.  The space of sections of
$\LL_\Sigma$ is certainly contractible, and hence simply-connected.
One can in
fact make sense of the space of $2$--cycles~\cite{almgren:cycles}, and its
fundamental group turns out to be precisely $H_3(Y)$.   For our
purposes, though, we do not need this remarkable fact, and will work directly
with the condition that $H_3(Y) = 0$.   By Poincar\'e duality, this
is equivalent to assuming $H^1(Y) = 0$.  For simplicity, we will in fact assume
that $\pi_1(Y) = 0$.  Thus $H_1(Y)=0$, which in turn is the condition needed to
incorporate the $0$--dimensional class. 

We will initially define a polynomial $D(f)$ of degree $d$, under one of two
hypotheses.  We  assume that either $w_2(P)\neq 0$ , or that invariants are
being computed in the `strong' stable range: $d \ge 2c_2(P) + 2$.   Either of
these assumptions will ensure, via the standard counting argument of Donaldson
theory, the compactness of all of the low-dimensional moduli spaces which appear
in the definition.  In Section~\ref{basic}, we will prove a blow-up formula,
which will then be used to define $D(f)$ in all degrees.

For a collection of homology classes $[\Sigma_i]\in H_2(Y)$, represented by
embedded surfaces $\Sigma_i$, one could consider divisors $V_{\Sigma_i}$, in
sufficient numbers so that the moduli space
$$
\left(\bigcup_t \MM (P;\{g_t\})\right) \cap  \left(\bigcap_i
V_{\Sigma_i}\right)
$$
is $0$--dimensional.  (The $0$--dimensional class could be included, in a similar
manner.)  Let us write
$$ D(P;\{g_t\},V_{\Sigma_1},\dots,V_{\Sigma_d})
$$
for the algebraic count of points in this intersection.
$ D(P;\{g_t\},V_{\Sigma_1},\dots,V_{\Sigma_d}) $ is readily seen to be
independent of the choice path connecting $g_0$ and $f^*g_0$, by the same
argument as outline above.   However, the argument that this count is
independent of representatives of the divisors and of initial metric
$g_0$ breaks down.   To see this (and what to do about it) consider, as in the
discussion above two initial metrics
$g_0=K_{0,0}$ and $k_0=K_{0,1}$, with a generic path $K_{0,t}$ between them.
Following that construction, we take a $2$--parameter family of
metrics $K_{s,t}$
(with $K_{1,t} = f^*K_{0,t}$)
and an associated $2$--parameter moduli space $
\Tilde{\Tilde{\mathcal{M}}}^{YM}$.  Intersecting with the divisors $V_{\Sigma}$
gives a null-cobordism of $\partial\Tilde{\Tilde{\mathcal{M}}}^{YM}$.
{\em A priori}, $f$ does not match up the right and left sides of this
cobordism, as one would need in order to get a cobordism between top and
bottom.

Indeed,  $f$ induces an isomorphism between
\begin{equation}\label{modiso}
\bigcup_t \MM(P;K_{0,t}) \cap  \left(\bigcap_i V_{\Sigma_i}\right)\quad
\text{and}
\quad
\bigcup_t \MM(P;K_{1,t}) \cap  \left(\bigcap_i f^* V_{\Sigma_i}\right)
\end{equation}
where $f^*V_\Sigma$ is the inverse image of $V_\Sigma$ under the diffeomorphism
$f^* \co \BB(P) \to \BB(P)$ induced by $f$.  There is no good reason to
expect that
$f^* V_{\Sigma} = V_\Sigma$.  Among other things,
$f(\Sigma)$ might not even be homologous to $\Sigma$.  In order to get a
diffeomorphism invariant, some restrictions are needed; here is one approach.

Let $\VV$ (or $\VV(f)$ if the diffeomorphism needs to be specified) be the
subgroup of
$H_2(Y)$ fixed by the action of
$f_*$; the invariant will be a polynomial in $H_0(Y) \oplus \VV$.  Represent an
element in $\VV$ by a generic surface $\Sigma$ in $Y$, and choose a
generic $3$--chain
$C$ giving a homology between $\Sigma$ and $f(\Sigma)$.  (From a technical
point of view, it would perhaps be preferable to let $C$ be the image of an
oriented
$3$--manifold via a smooth map to $Y$, but we will ignore this point for the
moment.)
As in Donaldson's original work~\cite{donaldson:polynomial}, consider a
line bundle
$\LL_\Sigma\to \BB^*(C)$ and a section $s_\Sigma$ whose
pull-back to $\BB(Y)$ defines the divisor $V_\Sigma$.  Using the action of $f$,
we get a section of $\LL_{f(\Sigma)}$, whose
divisor is $V_{f(\Sigma)}$.  Now $\LL_\Sigma$ and $\LL_{f(\Sigma)}$ are
equivalent when pulled back to $\BB(C)$, and a choice of homotopy between their
corresponding sections gives a cobordism $V_C$ between
$V_\Sigma$ and $V_{f(\Sigma)}$.

By analogy with $D(P;\{g_t\},V_{\Sigma_1},\dots,V_{\Sigma_d})$, we define, for
any $3$--chain $C_1$ and generic metric $g$, the invariant
$$
  D(P;g,V_{C_1},V_{\Sigma_2},\dots,V_{\Sigma_d}) =
\#\left[\MM(P;g)\bigcap \left( V_{C_1}\cap V_{\Sigma_2}\cdots \cap
V_{\Sigma_d}\right) \right].
$$
The term corresponding to the $3$--chain can go in any slot, in place of the
corresponding $V_\Sigma$.

Using the cobordisms $V_C$, we can finally give the actual
definition of a polynomial invariant.  
\begin{defn}\label{polydefinition}
Let $f\co Y \to Y$ be an orientation preserving diffeomorphism.  Assume that:
\begin{enumerate}
\item $H_1(Y) = 0$.
\item $P$ is an $\SO(3)$ or $\SU(2)$ bundle such that $f^*P \cong P$.
\item $w_2(P) \neq 0$ or (if $P$ is an $\SU(2)$ bundle) $d \ge 2c_2(P) +2$.
\item $\alpha(f) \beta(f) = 1$.
\end{enumerate}
Let $\Sigma_1,\ldots,\Sigma_d$ be generic surfaces carrying homology classes in
$\VV=
\ker(f_* -1)$, and suppose that $-2p_1(P) - 3(b_2^+(Y) + 1)  = 2d -1$.  Let
$C_1,\ldots,C_d$ be generic
$3$--chains in $Y$ such that $\partial C_i = f(\Sigma_i) - \Sigma_i$.  For
a metric $g_0$ on $Y$, let $\{g_t\} $ be a smooth path such that $g_1
= f^*g_0$.
Define
\begin{equation}\label{invdef}
\begin{split}
D_Y(f;\Sigma_1,\ldots,\Sigma_n)  = &
D(P;\{g_t\},V_{\Sigma_1},\dots,V_{\Sigma_d}) \\
& +
D(P;g_1,V_{C_1},V_{\Sigma_2},\dots,V_{\Sigma_d}) \\
& +
D(P;g_1,V_{f(\Sigma_1)},V_{C_2},\dots,V_{\Sigma_d}) \\
& \hspace*{.1in}\vdots\\
& +
D(P;g_1,V_{f(\Sigma_1)},V_{f(\Sigma_2)},\dots, V_{C_d}).
\end{split}
\end{equation}
\end{defn}
The term `generic' for a metric $g$ means that the moduli space is smooth of
the expected dimension, with no reducibles.   All surfaces $\Sigma_i$ and
$3$--chains $C_j$, as well as sections of associated line bundles (and
homotopies of such) are to be in general position, so the intersections with
$\MM(P;g)$ is smooth of the expected dimension as well.  Without loss of
generality, one can demand that the same is true of intersections with divisors
$V_{f(\Sigma)}$ and $V_{f(C)}$ as well.
\begin{remark}
This definition seems complicated, so some explanation may be helpful.    The
idea of the invariants under discussion is to use the moduli space associated
to a path in the space of choices of parameters used in defining an ordinary
invariant of a single $4$--manifold.  The parameter space involved in the usual
degree--$d$ Donaldson invariant is roughly $\met(Y) \times (\mathcal{C}_2)^d$
where $\mathcal{C}_2$ is the space of $2$--cycles in the relevant homology
classes.  The role of a path in the
$k^{th}$ factor of
$\mathcal{C}_2$ is played by a $3$--chain $C_k$.  In these terms, the definition
says to take a `path' from $(g_0,\Sigma_1,\ldots,\Sigma_d)$ to
$(f^*g_0,f(\Sigma_1),\ldots,f(\Sigma_d))$ which is a composition of paths, each
having non-constant projection into one factor at a time.
\end{remark}

\begin{theorem}\label{polydef}
Under hypotheses (1)--(4) in Definition~\ref{polydefinition},
$D_Y(f;\Sigma_1,\ldots,\Sigma_n)$ does not
depend on the choice of initial generic metric $g_0$ and path $g_t$, on the
choice of surfaces representing $[\Sigma_i]$, or on the choice of $3$--chains
$C_i$.
\end{theorem}
\begin{proof}
The independence of $D(f)$ from choice (relative to the endpoints) of
the path
$g_t$ is identical to that given before, because the only term which could
possibly change is the first.  The independence from the initial point $g_0$ is
more elaborate, as suggested by the discussion above.  Let $k_0$ be another
generic metric, and $K_{s,t}$ a $2$--parameter family of metrics with
\begin{itemize}
\item $K_{0,t}$ a generic path from $g_0$ to $k_0$;
\item $K_{s,0}$ = a generic path from $g_0$ to $g_1=f^*g_0$;
\item $K_{s,1}$ = a generic path from $k_0$ to $k_1=f^*k_0$;
\item $K_{1,t}= f^*K_{0,t}$.
\end{itemize}
As before, we get a $2$--parameter moduli space
\addtocounter{equation}{1}
\begin{equation}\label{mod0}
\Tilde{\Tilde{\MM}}(P;\{K_{s,t}\})  = \left(\bigcup_{(s,t)\in I \times I}
\MM(P;K_{s,t})\right)
\bigcap
\left(V_{\Sigma_1} \cap V_{\Sigma_2}\cdots \cap V_{\Sigma_d}\right)
\tag{\theequation.0}
\end{equation}
which is a compact oriented $1$--manifold.
\addtocounter{equation}{-1}

Treating the $3$--chains $C_j$ as parameters, in the spirit of the preceding
remarks, we consider the following collection of
$2$--parameter moduli spaces, which again are $1$--dimensional manifolds with
boundary.
\begin{subequations}
   \renewcommand{\theequation}{\theparentequation.\arabic{equation}}%
\begin{equation}
\begin{align}\label{mod1}
\Tilde{\Tilde{\MM}}(P;\{K_{1,t}\},C_1)& =  \tilde\MM(P;\{K_{1,t}\}) \bigcap
  \left( V_{C_1}\cap V_{\Sigma_2}\cdots \cap
V_{\Sigma_d}\right) \\
\Tilde{\Tilde{\MM}}(P;\{K_{1,t}\},C_2)& =  \tilde\MM(P;\{K_{1,t}\}) \bigcap
\left( V_{f(\Sigma_1)}\cap V_{C_2}\cdots \cap V_{\Sigma_d}\right) \\
\vdots\hspace*{.5in} \ &\hspace*{.1in}\vdots\hspace*{1.05in}{\vdots}
\tag*{\vdots\hbox to .15in{}}
\stepcounter{equation} \\
\label{modd}
\Tilde{\Tilde{\MM}}(P;\{K_{1,t}\},C_d)& =  \tilde\MM(P;\{K_{1,t}\}) \bigcap
\left( V_{f(\Sigma_1)}\cap V_{f(\Sigma_2)}\cdots
\cap V_{C_d}\right)\tag{\theparentequation.d}
\end{align}
\end{equation}
\end{subequations}

The boundary of each of the $1$--dimensional moduli
spaces~\eqref{mod0},~\eqref{mod1},\dots,~\eqref{modd} has algebraically $0$
points.  As discussed before, the boundary of each $2$--parameter moduli
space can alternatively be described as the sum of the algebraic counts of
points in appropriate $1$--parameter moduli spaces.  This leads to $d+1$
equations:
\begin{subequations}
   \renewcommand{\theequation}{\theparentequation.\arabic{equation}}%
\begin{equation}
\begin{align}
0 = &
\ D(P;\{k_s\},V_{\Sigma_1},\dots,V_{\Sigma_d}) -
D(P;\{g_s\},V_{\Sigma_1},\dots,V_{\Sigma_d})\tag{\theparentequation.0}
\\
& -D(P;\{K_{0,t}\},V_{\Sigma_1},\dots,V_{\Sigma_d})
+ D(P;\{K_{1,t}\},V_{\Sigma_1},\dots,V_{\Sigma_d})\notag
\\
0 = &
\ D(P;k_1,V_{C_1},V_{\Sigma_2},\dots,V_{\Sigma_d})
- D(P;g_1,V_{C_1},V_{\Sigma_2},\dots,V_{\Sigma_d})
\\
&
-D(P;\{K_{1,t}\},V_{\Sigma_1},V_{\Sigma_2},\dots,V_{\Sigma_d})
+ D(P;\{K_{1,t}\},V_{f(\Sigma_1)},V_{\Sigma_2},\dots,V_{\Sigma_d})\notag
\\
& \hspace*{1in}\vdots\tag*{\hspace*{-.05in}\vdots\hbox to .15in{}}
\\
0 = &
\ D(P;k_1,V_{f(\Sigma_1)},V_{f(\Sigma_2)},\dots, V_{C_d})
- D(P;g_1,V_{f(\Sigma_1)},V_{f(\Sigma_2)},\dots,V_{C_d})
\tag{\theparentequation.d}\\
&
\hspace*{-30pt}-D(P;\{K_{1,t}\},V_{f(\Sigma_1)},V_{f(\Sigma_2)},\dots,V_{\Sigma_d})
+ D(P;\{K_{1,t}\},V_{f(\Sigma_1)},V_{f(\Sigma_2)},\dots,V_{f(\Sigma_d)})\notag
\end{align}
\end{equation}
\end{subequations}
Adding these equations together, most of the terms cancel in pairs, leaving
the difference between the invariant calculated with the paths $\{k_s\}$ and
$\{g_s\}$, plus
$$
D(P;\{K_{1,t}\},V_{f(\Sigma_1)},V_{f(\Sigma_2)},\dots,V_{f(\Sigma_d)}) -
D(P;\{K_{0,t}\},V_{\Sigma_1},\dots,V_{\Sigma_d}).
$$
However, the isomorphism~\eqref{modiso}, coupled with the orientation
hypothesis that
$\alpha(f)\beta(f) =1$, means that the two terms are equal, and so
the invariant
doesn't depend on the choice of initial metric $g_0$.

The other choices of parameters involved in the definition of $D(f)$ are: the
specific surface representing $[\Sigma_i]$, the choice of section defining
$V_{\Sigma_i}$, the choice of $3$--chain $C_i$ with $\partial C_i =
f(\Sigma_i) -
\Sigma_i$, and the section defining $V_{C_i}$.  As remarked earlier, the
verification that, for fixed $\Sigma_i$, the choices of section don't
affect the
value of $D(f)$ is virtually identical to arguments given above, because
sections vary in a contractible space.  A similar remark applies to the
choice of $V_C$, given a specific $3$--chain $C$.

The independence from the choice of $\Sigma$'s and $C$'s differs in that a
substitute must be found for one basic mechanism: the existence of the family
$K_{s,t}$ derives from the fact that the space of metrics is simply
connected.  The idea is the same for all of the choices; we will illustrate the
point in the simplest instance.   So suppose that two $3$--chains $C_1$ and
$C_1'$ are given, both of which have boundary $f(\Sigma_1)-\Sigma_1$.  The only
place in equation~\eqref{invdef} in which $C_1$ enters is in the term
$$
D(P;g_1,V_{C_{1}},V_{\Sigma_2},\dots,V_{\Sigma_d}).
$$
Because the $3$--chains have the same boundary, it follows that $C_1' - C_1$
is a $3$--{\em cycle} which is a boundary of a $4$--chain $\Delta$, by our
hypothesis that $H_3(Y) =0$.  One can use restriction to connections on
$\Delta$ to define a $1$--dimensional moduli space
$\Tilde{\Tilde{\mathcal{M}}}$.  Taking the boundary of this moduli space gives
$$
D(P;g_1,V_{C_1},V_{\Sigma_2},\dots,V_{\Sigma_d}) =
D(P;g_1,V_{C_1'},V_{\Sigma_2},\dots,V_{\Sigma_d})
$$
by the standard argument.
\end{proof}

A similar technique may be used to incorporate the
$0$--dimensional class.  The invariant is readily checked to be multilinear, and
so defines a polynomial invariant in $P[H_0(Y) \oplus \VV(f)]$.  Some other
basic properties are summarized in the following theorem; they are analogous to
properties which hold for the degree $0$ part, and are proved in the same way.

\begin{theorem}\label{properties} Let $f$ and $g$ be diffeomorphisms for
which invariants
$D_Y(f)$ and $D_Y(g)$ are defined.
\begin{enumerate}
\item
The polynomials of a composition are defined on $H_0(Y) \oplus \VV(f,g)$, where
$\VV(f,g) = \VV(f) \cap \VV(g)$, and satisfy
$$
D_Y(f\circ g) = D_Y(g \circ f) = D_Y(f) + D_Y(g).
$$
\item The polynomial of $f^{-1}$ is $-D_Y(f)$.
\item If $f$ and $g$ are isotopic, then $D(f) = D(g)$.
\end{enumerate}
\end{theorem}

Because the applications are all to simply-connected manifolds, we haven't
stated the theorems in maximum generality.  The weakest set of hypotheses which
would give rise to an invariant of the type described in this section would
seem to be that $H_1(Y;\mathbf{Q}) =0$, and that $w_2(P)$ is not the pullback
of a class in $H^2(B\pi_1(Y);\mathbf{Z}_2)$.  The invariant would then be 
$\mathbf{Q}$ rather than in $\mathbf{Z}$--valued.

\section{Some basic theorems of $1$--parameter gauge theory}\label{basic}
In this section we will state (and sketch proofs of) analogues of the basic
connected-sum and blowup formulas for the Donaldson invariant.  Undoubtedly,
more elaborate versions of the gluing principles in gauge theory will work in
the $1$--parameter context, but we will state only those theorems which we
actually use.   The simple situation in which we work may be summarized in the
following definition.

\begin{defn}\label{diffsum}
Suppose that $f$ and $g$ are diffeomorphisms of manifolds $X$ and $Y$, which
are the identity near base points $x$ and $y$.  The connected sum $f\# g$ is
the obvious diffeomorphism on the connected sum $X\# Y$; it depends up to
isotopy only on the isotopy classes of $f$ and $g$ relative to neighborhoods of
the base points.
\end{defn}

A useful (cf~\cite{MMnote}) technical device for the ordinary Donaldson
polynomial is the fact that no (rational) information is lost if one replaces a
manifold by its connected sum with $\cpbar$.   A similar principle holds for
the $1$--parameter invariants.  To state this, let $L_0 \to \cpbar$ be the
complex line bundle such that $c_1(L_0)$ is Poincar\'e dual to the exceptional
curve
$E$ in $\cpbar$.
\begin{theorem}\label{blowup}
Suppose that the polynomial invariant $D(f,P)$ is defined for a diffeomorphism
$f\co Y \to Y$.  Then the invariant $D(f\# id_{\cpbar},P\#(L_0 \oplus \R))$ is
defined, and satisfies
\begin{equation}\label{blowupformula}
D(f\# id_{\cpbar},P\#(L_0 \oplus \R))(E,E) = -2 D(f,P).
\end{equation}
\end{theorem}
\begin{proof}
Choose a path of metrics and a collection of $3$--chains $C_i$ with $\partial
C_i = f(\Sigma_i) - \Sigma_i$ which define $D(f,P)$.  The path can assumed to
be constant near the connected sum point, so it extends to give a path of
metrics on $Y\# \cpbar$.  Similarly, the $3$--chains can be assumed to miss the
connected sum point, so they are $3$--chains in the connected sum in a natural
way.

Now we use a standard gluing argument:  choose a metric on $Y\# \cpbar$ with a
long tube along the $S^3$.  For sufficiently long tube length, we can calculate
each term in the definition of the invariant.   The $3$--chain  $C$ with
$\partial C= f(C) -C$ may be taken to be degenerate, so that the last two terms
(of the
form $D(P;g_1,V_{f(\Sigma_1)},V_{f(\Sigma_2)},\dots,V_{f(\Sigma_d)},V_{f(E)},
V_{C})$) are $0$ for dimensional reasons.  The moduli spaces corresponding to
the other terms in the definition, may all be described by the Kuranishi model
for the $1$--parameter moduli space, as in~\cite{ruberman:isotopy}.  The local
picture, and hence the calculation of the coefficients, is the same as in the
proof of the usual blowup formula.
\end{proof}

Following the scheme laid out in~\cite{MMnote}, we can extend the definition of
the invariants $D(f)$ outside the `stable range' by repeatedly blowing up to
increase the energy, and then using~\eqref{blowupformula}. The result is that
the invariant of $f$ is a collection of rational-valued polynomials
of arbitrary
degree in $H_0(Y) \oplus \VV(f)$.
Following~\cite{kronheimer-mrowka:relations}
we introduce the notion of a diffeomorphism being of {\em simple type}, and
assemble the polynomials into a formal power series $\DD(f)$, which we will
call the Donaldson series of $f$.  In the examples discussed below and in the
next section, the power series are determined by a set of basic classes, as in
the main theorem of~\cite{kronheimer-mrowka:relations}.  It would be
of interest
to know if such a structure theorem holds more generally under the simple type
hypothesis.

There is also a version of the connected sum theorem; the proof is a
simple dimension-counting argument and will be omitted.
\begin{theorem}\label{sum}
Suppose that $f_i\co Y_i \to Y_i$ are diffeomorphisms, where $Y_i$\break (for
$i=1,2$) are
$4$--manifolds satisfying $b^2_+(Y_i) \ge 2$.
Then any invariant\break
 $D(f_1\# f_2,P_1\# P_2)$ which is defined must vanish.
\end{theorem}

The remaining case to investigate is when $b^2_+(Y_1) \ge 2$ and $b^2_+(Y_2)
=1$.  The result is more complicated, and it depends on the behavior of the
diffeomorphism $f_2$.  The basic idea is that the evaluation of the
$1$--parameter invariant on homology classes supported in $Y_1$ is, in some
circumstances, the product of an ordinary Donaldson invariant of $Y_1$ with a
term related to the wall-crossing  phenomenon characteristic of gauge theory on
manifolds with
$b^2_+ = 1$.  A completely general treatment would run into unresolved problems
associated with that theory (under the general rubric of the {\em
Kotschick--Morgan conjecture}---cf~\cite{kotschick:so3,kotschick-morgan:so3}).
We will state a relatively simple version, which avoids these technicalities,
but which suffices for the main application.  A reasonable extension of this
statement, parallel to the Kotschick--Morgan conjecture, would be that the
restriction of
$D(f_1\#f_2)$ to $H_2(Y_1)$ depends in some universal fashion on $D(f_1)$ and
the action of $f_2$ on cohomology.  The full polynomial (ie including
$H_2(Y_2)$) is also of interest.

Let $N$ be a simply-connected manifold with $b_+^2 = 1$, and let $L \to N$ be a
complex line bundle with $c_1(L)^2 = -1$.  Note that this implies
that $w_2(P_N)
\neq 0$, where $P_N$ is the $\SO(3)$ bundle over $N$ associated to $L
\oplus \R$.  A choice of orientation for $H^2_+(N)$ picks out a positive sheet
of the hyperboloid $\HH =\{ \alpha \in H^2_+(N) \, | \, \alpha^2 =1
\}$.  Inside
$\HH$ lie the walls $\WW$, where a wall is the orthogonal complement
(intersected with $\HH$) of a class $x \in H^2(N;\Z)$ satisfying $x\equiv
c_1(L) \pmod 2$ and $x^2 = 1$.   The walls are transversally
oriented, and form a
locally finite something or other.   Note that any metric $g$ on
$N$ determines a unique self-dual harmonic form $ \omega_g \in \HH$, called its
period point.

Let $f_N$ be a diffeomorphism of $N$ which is the identity near a point of
$N$, which has the property that
$f_N^*$ preserves $w_2(P_N)$, and satisfies $\alpha(f_N)\beta(f_N) = 1$. Such
diffeomorphisms were constructed on $N = \cp^2 \#_2\cpbar$ in \S3
of~\cite{ruberman:isotopy}, and easily extend to arbitrary connected sums
$\cp^2 \#_k\cpbar$.  Let $g^N_0$ be a metric on $N$, which is fixed by $f$ near
the connected sum point, and whose period point does not lie on any of the
walls.  Join $g^N_0$ to $f^*(g^N_0)$ by a path whose induced path $\gamma$ of
period points is transverse to $\WW$.  Using the transverse orientation of
$\WW$, the intersection number of this path with $\WW$ is well-defined.
\begin{theorem}\label{gluing}
Let $f$ be the diffeomorphism of $Z= Y \# N$ gotten by gluing $f_N$ to the
identity of $Y$.  Then $D_Z(f) = 2(\gamma\cdot\WW) D_Y$.
\end{theorem}

\section{Applications to the topology of the diffeomorphism group}

The $1$--parameter invariants, as extended in the previous section, fit together
naturally to give a homomorphism which will show that $\pi_0(\diffh)$ can be
infinitely generated, proving Theorem~\ref{infgen} of the
introduction.  (Recall
that $\diffh(Z)$ is the subgroup of the diffeomorphism group consisting of
diffeomorphisms which act trivially on homology.)

There is a small technical observation to be made in order to draw conclusions
about $\pi_0(\diff)$ from  our results.  Namely, two diffeomorphisms are
in the same path component of $\diff$ if and only if they are isotopic.  This
seems a little surprising at first, because there is no smoothness required for
a path in $\diff$.  The proof relies on simple properties of the Whitney
$\mathcal{C}^\infty$ topology on smooth maps, and is quite standard in the
subject---compare~\cite[Definition 3.9 and Problem 4.6]{munkres:diff}
and~\cite{burghelea-lashof-rothenberg}.

Combining this observation, the definition of the Donaldson series of a
diffeomorphism, and Theorem~\ref{properties}, we get the following result.
\begin{theorem}\label{power series}
Let $Y$ be a $4$--manifold with $b_+^2$ an even number $\ge 4$.  Then the
Donaldson series defines a homomorphism
$$
\D\co  \pi_0(\diffh(Y)) \rightarrow \R[[ H_2(Y)^* ]].
$$
\end{theorem}

The proof of Theorem~\ref{infgen} will be completed by showing that for the
manifolds
$Z_n$ described in the introduction, the image is infinitely generated.
\begin{proof}[Proof of Theorem~\ref{infgen}]
Suppose that $Z$ is of the form $Y\# N$, where $N=\cp^2\#_2\cpbar$, and notice
that restriction defines a homomorphism 
$$
r_Y^*\co  \R[[ H_2(Z)^* ]] \rightarrow \R[[ H_2(Y)^* ]].
$$
Let $f$ be a diffeomorphism of the form $id_Z\#f_N$,  as discussed before
Theorem~\ref{gluing}.   In particular, $f_N$ should be chosen as a composition
of reflections in two different $(-1)$--spheres, as in~\cite{ruberman:isotopy};
the intersection number $\gamma\cdot \WW$ is computed in that paper to be $-2$.
Suppose finally that
$Y$ has  simple type in the sense
of~\cite{kronheimer-mrowka:relations}, so that
its Donaldson series
$\D_Y$ is determined by a finite set of  basic classes $\kappa_i(Y) \in
H^2(Y;\Z)$.  Rewriting Theorem~\ref{gluing} in terms of the Donaldson series of
$f$, we see that $r_Y^*\D(f) = -4\D_Y$.
In particular,
$r_Y^*\D(f)$ has the form described by the structure theorem
of~\cite{kronheimer-mrowka:relations}, and so is determined by the same set of
basic classes $\kappa_i(f) =
\kappa_i(Y)$.  Moreover, the coefficients $\beta_i(f)$ in the expression of the
series as a sum of exponentials of the $\kappa_i$, are equal to the
corresponding coefficients for $Y$.

Under composition of diffeomorphisms the Donaldson series
add.  For diffeomorphisms $f,g \in \diffh(Z)$ whose series are determined by
basic classes,  this implies the following statement.  The set of basic classes
for $f \circ g$ is the union of the set of basic classes for $f$ and for $g$,
leaving out those basic classes which $f$ and $g$ have in common but whose
coefficients cancel.  In other words, a basic class $\kappa_i(f) = \kappa_j(g)$
is removed from the union if the coefficient
$\beta_i(f) = -\beta_j(g)$.

In the paragraphs which follow, we will show that if $Z$ is any one of the
manifolds described in the statement of Theorem~\ref{infgen}, then it admits a
series of diffeomorphisms $\{f_j \}\ (j=1,\dots,\infty)$ which are
all homotopic
to the identity, with the property that $f_m$ has at least $m$ different basic
classes.  We claim that the image under $\D$ of the subgroup of
$\pi_0(\diffh(Z))$ generated by the $f_j$ is infinitely generated.  Suppose
that the diffeomorphisms have been indexed so that $f_m$ has at least one basic
class which does not occur in the list of basic classes for the $f_j$ for $j <
m$.  Note that if
$K_1,\ldots,K_n$ are distinct elements in $H^2(Y)$, then the exponentials
$\exp(K_1),,\ldots,\exp(K_n)$ are linearly independent elements in the power
series ring $\R[[ H_2(Y)^* ]]$.   Thus in any any linear relation
$$
\sum_{j=1}^m a_j\D(f_j) = 0
$$
the coefficient $a_m $ must be $0$.  The claim follows immediately by
induction, and so we have that $\pi_0(\diffh(Z))$ is infinitely generated.

Let $Y_n$ denote $\#_{2n-1} \cp^2 \#_{10n-1} \cpbar$ for $n$ odd, and
$\#_{2n-1}
\cp^2 \#_{10n} \cpbar$ for $n$ even.  The manifold $Z_n$ will be
simply $Y_n\# N$,
where $N = \cp^2\#_2\cpbar$ as before.  Let $E(n)$ be the
simply-connected elliptic
surface with $p_g = n-1$ and no multiple fiber, and let $E(n;p)$ denote the
result of a single logarithmic transform on a fiber in $E(n)$.  The standard
convention is that $E(n;1)$ is the same as $E(n)$.

We will make use of the following facts about these manifolds.
\begin{enumerate}
\item For $n$ odd, $E(n;p)\simeq Y_n$, and for $n$ even, $E(n;p)\#
\cpbar \simeq
Y_n$.
\item   $E(n;p)\# \cp^2$ decomposes completely into a connected sum
of $\cp^2$'s
and $\cpbar$'s.  See~\cite{mandelbaum:survey}
or~\cite{mandelbaum:elliptic,moishezon:sums} for more details.
\item The diffeomorphism group of $Y_n\# \cp^2$ acts transitively on elements
in $H_2(Z_n)$ of given square, divisibility, and type (ie characteristic or
not)~\cite{wall:diffeomorphisms}.
\item The Donaldson series for $E(n;p)$ is given~\cite{fs:rational} by
$$
\D_{E(n;p)} = \exp(Q/2)\frac{\sinh^{n-1}(f)}{\sinh(f_p)}
$$
where $f_p$ is the multiple fiber (and therefore the regular fiber $f=pf_p$ in
homology).
\item The Donaldson series for  $E(n;p)\#\cpbar$ is $\D_{E(n;p)}
e^{-\frac{E^2}{2}}\cosh(E)$ where $E$ is dual to the exceptional class.
\end{enumerate}

The argument differs in minor details between the cases when $n$ is even or
odd; for simplicity we will concentrate on $n$ odd.  The main point of this is
that $E(n;p)$ is not spin when $n$ is odd.

Let $S_0$ denote the standard  (complex) $2$--sphere in
$\cp^2$, viewed as a submanifold in $Y_n\# \cp^2$, and let $S_p'$ denote the
analogous sphere in $E(n;p)\# \cp^2$.  Using the first two items, choose a
diffeomorphism of $E(n;p)\# \cp^2$ with $Y_n\# \cp^2$.   Since $S_p'$
is not characteristic, any initial choice of diffeomorphism may be varied
by a self-diffeomorphism of $Y_n\# \cp^2$ to ensure that the image of $S_p'$ is
homologous to $S_0$.  Denote this image, viewed as a sphere in $Y_n\# \cp^2$ or
in $Z_n$, by $S_p$.    Note that the homology of $Y_n$ may be identified with
the orthogonal complement to $S_p$, with respect to the intersection pairing,
and hence the image of $H_2(E(n;p))$ is $H_2(Y_n)$.

As in~\cite{ruberman:isotopy}, the $(-1)$--spheres $S_p \pm E_1 + E_2$ in $Z_n$
give rise to reflections $\rho_p^\pm$, and we set
$$
f_p = (\rho_p^+ \circ \rho_p^-) \circ  (\rho_0^+ \circ \rho_0^-)^{-1}.
$$
Because $S_p$ and $S_0$ are homologous, the action of $f_p$ on homology is
trivial, and thus~\cite{quinn:isotopy,cochran-habegger:homotopy} $f_p$ is
homotopic to the identity.  The image of $\D(f_p)$ under $r_Y^*$ is the
Donaldson series of $E(n;p)$, and so is given by the formula in item 4 above.
Expanding the hyperbolic sines, we see that  $E(n;p)$ has $(n-1)p$ basic
classes, and so there are the same number of basic classes for
$r_Y^*(\D(f_p)$.  Thus the $f_p$ generate an infinitely generated subgroup of
$\diffh(Z_n)$.
\end{proof}

\Addresses\recd

\end{document}